\documentclass{llncs}
\usepackage{algorithm}
\usepackage{algorithmic}

\usepackage{amsmath}
\usepackage{amsfonts}
\usepackage{epsfig}
\usepackage{graphics}

\newcommand{\reals}{\mathbb{R}}

\begin{document}

\mainmatter

\title{A hybrid constraint programming and semidefinite
  programming approach for the stable set problem}


\author{W.J. van Hoeve\inst{}}

\institute{CWI, P.O. Box 94079, 1090 GB Amsterdam, The Netherlands\\
\email{w.j.van.hoeve@cwi.nl}\\
\texttt{http://www.cwi.nl/\~{ }wjvh/}}

\maketitle
\begin{abstract}
This work presents a hybrid approach to solve the maximum stable set
problem, using constraint and semidefinite programming. The approach
consists of two steps: subproblem generation and subproblem solution.
First we rank the variable domain values, based on the solution of a
semidefinite relaxation. Using this ranking, we generate the most
promising subproblems first, by exploring a search tree using a
limited discrepancy strategy. Then the subproblems are being solved
using a constraint programming solver. To strengthen the semidefinite
relaxation, we propose to infer additional constraints from the
discrepancy structure. Computational results show that the
semidefinite relaxation is very informative, since solutions of good
quality are found in the first subproblems, or optimality is proven
immediately.
\end{abstract}

\section{Introduction}
This paper describes a hybrid method to solve the maximum weighted
stable set problem, or {\em stable set problem}\footnote{ Alternative
  names for the stable set problem are vertex packing,
  coclique or independent set problem.}  in short. 
Given a graph with weighted vertices, the stable set
problem is to find a subset of vertices of maximum weight, such that
no two vertices in this subset are joined by an edge of the graph.
In the unweighted case (when all weights are equal to 1), this 
problem amounts to the maximum cardinality stable set problem, which
has been shown to be already NP-hard~\cite{papadimitriou_steiglitz}.

We propose a two-phase approach to solve the stable set problem,
either with or without proving optimality. The first phase generates
subproblems based upon a semidefinite relaxation, the second phase
solves the subproblems using constraint programming. Concerning the
first phase, given a model for the stable set problem, we solve its
semidefinite relaxation. The solution provides us fractional values
for the variables of the model. These fractional values are a good
indication for the optimal (discrete) values of the variables. Hence
we divide selected variable domains in two parts: a `good' subdomain
and a `bad' subdomain. By branching on these subdomains using a limited
discrepancy strategy~\cite{harvey95}, we obtain first a very promising
subproblem, and subsequently less promising subproblems.

The second phase consists of the solution of the subproblems. Since
they are much smaller than the original problem, we can easily solve
them using a constraint programming solver. 

As computational results will show, the semidefinite relaxation is
quite informative. In several cases we can simply round the solution
of the relaxation and obtain a provable optimal solution
already. Otherwise, we are likely to find a good solution in one of
the first subproblems. Using a limited number of subproblems to
investigate, we yield an {\em incomplete} method producing good
solutions. In order to obtain a {\em complete} search strategy,
we need, in principle, to generate and solve all possible
subproblems. A good upper bound is necessary to prove optimality
earlier. For this reason we investigated the use of `discrepancy
cuts' that can be added to the semidefinite program to strengthen the
relaxation and thus prune large parts of the search tree. However,
computational results will show that they can not be applied
efficiently on the instances we considered.

The outline of the paper is as follows. The next section gives a
motivation for the approach proposed in this work. Then, in
Section~\ref{sc:sdp} some preliminaries on semidefinite programming
are given. In Section~\ref{sc:formulations} we introduce the stable
set problem, integer optimization formulations and a semidefinite
relaxation. A description of our solution framework is given in
Section~\ref{sc:framework}. Section~\ref{sc:results} presents the
computational results. This is followed by an overview of related
literature in Section~\ref{sc:literature}. Finally, in
Section~\ref{sc:conclusions} we conclude and discuss future
directions. 

\section{Motivation}
Let us first motivate why a semidefinite relaxation is used rather
than a linear relaxation. Indeed, one could argue that linear programs
are being solved much faster in general. However, for the stable set
problem, linear relaxations are not very tight. Therefore one has to
identify and add inequalities that strengthen the relaxation. But it
is time consuming to identify such inequalities, and by enlarging the
model the solution process may slow down.

Several papers on approximation theory following~\cite{GW_maxcut} have
shown the tightness of semidefinite relaxations. However, being
tighter, semidefinite relaxations are more time-consuming to solve
than linear relaxations in practice. Hence one has to trade strength
for computation time. For 
large scale applications, semidefinite relaxations can often be
preferred as the relaxation of choice to be used in a branch and bound
framework. Moreover, our intention is not to solve a relaxation at
every node of the search tree. Typically, we only solve a relaxation
at the root node of the first phase (its solution is used to identify
the subdomains), and optionally at the root node of each subproblem
(in order to strengthen the upper bound). Therefore, we are willing to
make the trade-off in favour of the semidefinite relaxation.

Another point of view is the following. Although semidefinite
programming has been developing for many years now in the operations
research community, no efforts of integration or cooperation with
constraint programming have been made to our knowledge. Application
of semidefinite programming to problems typical to constraint
programming, as was done in the papers on approximation algorithms
mentioned in Section~\ref{sc:literature}, is not yet hybrid problem 
solving. In this paper a first step is being made. The solution of the
semidefinite relaxation is used to identify promising subdomains, and
also produces a tight upper bound for the constraint programming
solver. On the other hand, the solutions found by the constraint
programming solver serve as a lower bound inside the semidefinite
programming solver.

\section{Preliminaries on semidefinite programming}\label{sc:sdp}
Combinatorial optimization problems that are NP-hard are often solved
with the use of a relaxation. Typically, the integral restriction on
the decision variables is relaxed to obtain a continuous relaxation,
which can be polynomially solvable. Here we introduce semidefinite
programming as an extension of the more common linear
programming. Both paradigms can be used as a continuous relaxation.

In linear programming, combinatorial optimization problems are modeled
in the following way: 
\begin{equation}
\begin{array}{rl}
{\rm max}  & c^{\sf T}x \\
{\rm s.t.} & Ax \leq b \\
           & x \geq 0
\end{array}
\end{equation}
Here $x \in \reals^n$ is an $n$-dimensional vector of decision
variables and $c \in \reals^n$ a cost vector of dimension $n$. The ($m
\times n$) matrix 
$A \in \reals^{m \times n}$ and the $m$-dimensional vector $b \in
\reals^m$ define $m$ linear constraints on $x$. In other words, this
approach models problems using {\em nonnegative vectors} of variables. 

Semidefinite programming makes use of {\em positive semidefinite
  matrices} of variables instead of nonnegative vectors. A matrix $X
\in \reals^{n \times n}$ is said to be positive semidefinite (denoted
by $X \succeq 0$) when $y^{\sf T}Xy \geq 0$ for all vectors $y \in
\reals^n$. Semidefinite programs have the form 
\begin{equation}\label{eq:sdp_general}
\begin{array}{rll}
{\rm max}  & {\rm tr} (CX) \\
{\rm s.t.} & {\rm tr} (A_j X) = b_j & (j=1,\dots,m)\\
           & X \succeq 0.
\end{array}
\end{equation}
Here ${\rm tr}(X)$ denotes the {\em trace} of $X$, which is the sum of
its diagonal elements, i.e. ${\rm tr}(X) = \sum_{i=1}^{n} X_{ii}$. The
cost matrix $C \in \reals^{n \times n}$ and the constraint matrices
$A_j \in \reals^{n \times n}$ are supposed to be symmetric. The $m$ reals
$b_j$ and the $m$ matrices $A_j$ define again $m$ constraints. 

We can view semidefinite programming as an extension of linear
programming. Namely, when the matrices $C$ and $A_j \;(j=1, \dots, m)$
are all supposed to be diagonal matrices, the resulting semidefinite
program is equal to a linear program. In particular, then a
semidefinite programming constraint $A_j X = b_j$ corresponds to the
linear programming constraint defined by the $j$-th row of $A$ and the
$j$-th entry of $b$.

Applied as a continuous relaxation, semidefinite programming in
general produces solutions that are much closer to the integral
optimum than linear programming. Untuitively, this can be explained as
follows. The problem is `lifted' to a higher dimension, where one can
exploit more of the structure of the problem. Unfortunately, it is
not a trivial task to catch this structure, and otain a good
semidefinite program for a given problem. 

Theoretically, semidefinite programs have been proved to be
polynomially solvable using the so-called ellipsoid method (see for
instance~\cite{GLS88}). In practice, nowadays fast `interior point'
methods are being used for this purpose (see~\cite{alizadeh95} for an
overview). Being a special case of semidefinite programming, linear
programs are also polynomially solvable using an interior point
method. However, they are often solved with a special linear
programming solver, the Simplex method. Although this method can have
an exponential behaviour in theory, in practice it is the fastest
method known for solving linear programs. 

\section{The stable set problem} \label{sc:formulations}
In this section, the stable set problem is formally defined, and
formulated in two different ways as an integer optimization
problem. From this, a semidefinite programming relaxation is inferred.

\subsection{Definition}
Consider an undirected weighted graph $G=(V,E)$, where $V$ is the set
of vertices and $E$ the set of edges $\{(i,j)| i,j \in V, i \neq j \}$ of
$G$. Let $|V|=n$ and $|E|=m$. To each vertex $i \in V$ a weight $w_i
\in \reals$ is assigned. A {\em stable set} is a set $S \subseteq V$
such that no two vertices in $S$ are joined by an edge in $E$. The
{\em stable set problem} is the problem of finding a stable set of 
maximum total weight in $G$. This value is called the {\em stable set
  number} of $G$ and is denoted by $\alpha(G)$. The maximum
cardinality (or unweighted) stable set problem can be obtained by
taking all weights equal to 1. 


\subsection{Integer optimization formulation}
Let us first consider an integer linear programming formulation.
We introduce binary variables to indicate whether or not a vertex
belongs to the stable set $S$. So, for $n$ vertices, we have $n$
integer variables $x_i$ indexed by $i \in V$, with initial
domains $\{0,1\}$. In this way, $x_i=1$ if vertex $i$ is in the stable
set $S$, and $x_i = 0$ otherwise. We can now state the objective
function, being the sum of the weights of vertices that are in the
stable set $S$, as $\sum_{i=1}^n w_i x_i$. Finally, we define the
constraints that restrict two adjacent vertices to be both inside $S$
as $x_i + x_j \leq 1$, for all edges $(i,j) \in E$. Hence the
integer linear programming model becomes:
\begin{equation} \label{eq:ilp_form}
\begin{array}{rrrl}
\alpha(G) = & {\rm max} & \sum_{i=1}^{n} w_i x_i \\
&  {\rm s.t.} & x_i + x_j \leq 1 & \forall (i,j) \in E \\
&   & x_i \in \{0,1\}  & \forall i \in V.
\end{array}
\end{equation}

Another way of describing the same solution set is by the following
integer quadratic program
\begin{equation}\label{eq:quadratic}
\begin{array}{rrrl}
\alpha(G) = & {\rm max} & \sum_{i=1}^{n} w_i x_i \\
 & {\rm s.t.}& x_i x_j = 0 & \forall (i,j) \in E \\
 & & x_i^2 = x_i & \forall i \in V.
\end{array}
\end{equation}
Note that here the constraint $x_i \in \{0,1\}$ is replaced by $x_i^2
= x_i$. This quadratic formulation will be used below to infer a
semidefinite programming relaxation of the stable set problem.

In fact, both model (\ref{eq:ilp_form}) and (\ref{eq:quadratic}) can
be used as a constraint programming model. We have chosen the first
model, since the quadratic constraints take more time to propagate
than the linear constraints, while having the same pruning power.

\subsection{Semidefinite programming relaxation}
The integer quadratic program (\ref{eq:quadratic}) gives rise to a
semidefinite relaxation introduced by Lov\'asz \cite{lovasz79}, see
Gr\"otschel et al.~\cite{GLS88} for a comprehensive treatment. The
value of the objective function of this relaxation has been named the
{\em theta number} of a graph $G$, indicated by $\vartheta(G)$.
Let us start again from model (\ref{eq:quadratic}). As was indicated
in Section~\ref{sc:sdp}, we want to transform the current model that
uses a nonzero vector into a model that uses a positive semidefinite
matrix to represent our variables. In the current case, we can
construct a matrix $X \in \reals^{n \times n}$ by defining $X_{ij} =
x_i x_j$. Let us also construct a $(n \times n)$ cost matrix $W$ with
$W_{ii} = w_i$ for $i \in V$ and $W_{ij} = 0$ for all $i \neq
j$. Since $X_{ii} = x_i^2 = x_i$, the objective function becomes
tr$(WX)$. The edge constraints are easily transformed as $x_i x_j = 0
\Leftrightarrow X_{ij} = 0$. The first step in the 
transformation of model (\ref{eq:quadratic}) can now be made:
\begin{equation}\label{eq:step1}
\begin{array}{rrl}
 {\rm max} & {\rm tr}(WX) \\
 {\rm s.t.}& X_{ij} = 0 & \forall (i,j) \in E \\
 & x_i^2 = x_i & \forall i \in V.
\end{array}
\end{equation}
This model is still a quadratic program, although reformulated.
The problem remains how to model the last, very important,
constraint. We need a mapping of the diagonal entries $X_{ii} = x_i^2$
with the vector entries $x_i$. For this, we extend $X$ with another 
row and column that contain vector $x$ (both indexed by 0). For this
reason, define the $(n+1) \times (n+1)$ matrix $Y$ as 
\begin{displaymath}
Y = \left(
\begin{array}{cccc}
1 & x_1 & \cdots & x_n \\
x_1 \\
\hdots &  & X \\
x_n  \\
\end{array}\right)
\end{displaymath}
where the 1 in the leftmost corner of $Y$ is needed to obtain positive
semidefiniteness. In this case we can express the required mapping as
$Y_{ii} = \frac{1}{2} Y_{i0} + \frac{1}{2} Y_{0i}$ (note that $X$ and
$Y$ are symmetric), since then $x_i^2 = Y_{ii} = \frac{1}{2} Y_{i0} +
\frac{1}{2} Y_{0i} = x_i$. The final step in the transformation
consists of replacing the constraints on $X$ by constraints on $Y$. In
particular, instead of demanding $X$ to be a product of nonnegative
vectors, we restrict $Y$ to be a positive semidefinite matrix. Namely,
if the vector $x$ represents a stable set, then the matrix $Y$ is
positive semidefinite. However, not all positive semidefinite $Y$
matrices represent a stable set, in particular its values can take
fractional values. 

In order to maintain equal dimension to $Y$, a row
and a column indexed 0 should be added to $W$, all entries of which
containing value 0. Denote the resulting matrix by $\tilde{W}$. The
theta number of a graph $G$ can now be described as 
\begin{equation} \label{eq:theta1}
\begin{array}{rrcl}
\vartheta(G) = & {\rm max} & {\rm tr(\tilde{W}Y)} \\
& {\rm s.t.}& Y_{ii} = \frac{1}{2}Y_{i0} + \frac{1}{2}Y_{0i} & (i =
1,\dots,n) \\ 
& & Y_{ij} = 0 & \forall (i,j) \in E \\
& & Y \succeq 0.
\end{array}
\end{equation}
By construction, the diagonal value $Y_{ii}$ serves as an indication
for the value of variable $x_i$ $(i \in \{1, \dots, n\})$ in a maximum
stable set. In particular, this program is a relaxation for the stable
set problem, i.e. $\vartheta(G) \geq \alpha(G)$. Note that
program~(\ref{eq:theta1}) can easily be rewritten into the 
general form of program~(\ref{eq:sdp_general}). Namely, $Y_{ii} =
\frac{1}{2}Y_{i0} + \frac{1}{2}Y_{0i}$ is equal to tr$(AY)$ where the
$(n+1) \times (n+1)$ matrix $A$ consists of all zeroes, except for
$A_{ii} = 1$, $A_{i0} = -\frac{1}{2}$ and $A_{0i} = -\frac{1}{2}$,
which makes the corresponding $b$ entry equal to 0. Similarly for the 
edge constraints.

The theta number also arises from other formulations, different from 
the above, see \cite{GLS88}. In our implementation we have used the
formulation that has been shown to be computationally most efficient
among those alternatives \cite{gruber_rendl}. Let us introduce that
particular formulation. Again, let $x \in \{0,1\}^n$ be a vector of
binary variables representing a stable set. Define the $n \times n$
matrix $X$ as $X = \frac{1}{x^{\sf T} x} x x^{\sf T}$. Furthermore,
let the $n\times n$ cost matrix $W$ be defined as 
\begin{displaymath}
W_{ij} = \left\{
\begin{array}{ll}
w_i &  i = j \\
\sqrt{w_i  w_j} & i \neq j
\end{array}\right.
\end{displaymath}
for $i,j \in V$. Then the following semidefinite program
\begin{equation} \label{eq:theta2}
\begin{array}{rrrl}
\vartheta(G) = & {\rm max} & {\rm tr}(WX) \\
& {\rm s.t.}& {\rm tr} (X) = 1 \\
& & X_{ij} = 0 & \forall (i,j) \in E \\
& & X \succeq 0
\end{array}
\end{equation}
gives exactly the theta number of a graph. Again, it is not difficult
to rewrite this program into the general form of
program~(\ref{eq:sdp_general}). 

Program~(\ref{eq:theta2}) uses matrices of dimension $n$ and $m+1$
constraints, while program~(\ref{eq:theta1}) uses matrices of
dimension $n+1$ and $m+n$ constraints. This gives an indication why
program~(\ref{eq:theta2}) is computationally more efficient.

\section{Solution framework} \label{sc:framework}
\subsection{Overview}
The two-phase solution approach proposed here is similar to the one
described in \cite{milano_hoeve02}. In the first phase subproblems are
generated, which are being solved in the second phase. A general
overview of the method is presented in Algorithm~\ref{alg:dbs}.

\begin{algorithm}[t]
\caption{Solution framework}
\label{alg:dbs}
\begin{algorithmic}
\STATE{read problem}
\STATE{set maximum discrepancy}
\STATE{solve semidefinite program (\ref{eq:theta2})} $\rightarrow$
upper bound
\STATE{round solution of (\ref{eq:theta2})} $\rightarrow$ lower bound
\FOR{$i \in V_0$}
  \STATE{define $D_i^{\rm good}$ and $D_i^{\rm bad}$ using solution of
  (\ref{eq:theta2})}
\ENDFOR
\STATE{set discrepancy = 0}
\WHILE{lower bound $<$ upper bound {\bf and} discrepancy $\leq$ maximum 
discrepancy}
\STATE{generate subproblem using LDS branching strategy on $D_i^{\rm good}$
  and $D_i^{\rm bad}$}
\STATE{solve subproblem} $\rightarrow$ lower bound
\STATE{discrepancy = discrepancy + 1}
\ENDWHILE
\end{algorithmic}
\end{algorithm}
Let us first explain how we use the solution of the semidefinite
program~(\ref{eq:theta2}) to partition the domain $D_i$ of a variable 
$x_i$ into $D_i^{\rm good}$ and $D_i^{\rm bad}$ (for $i \in V$). As
was stated before, the solution of program~(\ref{eq:theta2}) assigns
fractional values between 0 and 1 to the variables. Naturally, if for
a variable $x_i$ the corresponding fractional value $Y_{ii}$ is close
to 1, we regard 1 to be a good value for variable $x_i$. More
specifically, we select the variable $x_i$ with the highest fractional
value $Y_{ii}$, set $D_i^{\rm good} = \{1\}$ and $D_i^{\rm bad} =
\{0\}$, and mark it as handled. Then we mark all its neighbours $j$
(with $(i,j) \in E$) as being handled, keeping their original
domain $D_j = \{0,1\}$. This procedure is repeated until all variables
are handled. For convenience, we partition $V$ into two
distinct sets $V_0 = \{i \in V \; | D_i^{\rm good} = \{1\} \}$ and
$V_1 = V \setminus V_0$. Here $V_1$ represents the set of neighbours
of $V_0$, with $D_j = \{0,1\}$.

In a similar way, we can use the solution of the semidefinite
relaxation to compute a first feasible integer solution. Namely,
follow the same procedure, but now instantiate the selected variable
$x_i = 1$ and set its neighbours $x_j = 0$. The
objective value of this feasible integer solution is in many cases
already equal to the (downward rounded) solution of the semidefinite
relaxation. In that case, we have found an optimal solution and
finish. In other cases, we can still use this first solution as a
lower bound to be applied during the solution of the subproblem.

Next, we explain how to generate subproblems using these subdomains.
The generation of subproblems makes use of a tree structure of depth
$|V_0|$ in which we branch on $D_i^{\rm good}$ versus  $D_i^{\rm
  bad}$. The tree is traversed using a limited discrepancy strategy
(LDS) \cite{harvey95}. LDS visits the nodes of a search tree
differently from depth-first search. It tries to follow a given
suggestion as good as it can. Branches opposite to the suggestion are
regarded as discrepancies and are gradually allowed to be
traversed. The first `run' of LDS doesn't allow any discrepancies, the
second allows one, and so on. This means that for a particular
discrepancy $k$, a path from the root to a leaf is allowed to consist
of maximally $k$ right branches.  Typically this method is applied
until a limited number of discrepancies is reached (say 2 or 3), which
yields an incomplete search strategy. In order to be complete, one has
to visit all nodes, up to discrepancy $d$ for a binary tree of depth
$d$.

In our case, the suggestion that should be followed are branches of
the kind $D_i^{\rm good}$, while branches $D_i^{\rm bad}$ are regarded
as discrepancies. Hence, our first subproblem is the subproblem defined by 
program~(\ref{eq:ilp_form}), and with $x_i \in D_i^{\rm good}$ for all
$i \in V_0$. The next $|V_0|$ subproblems have all $x_i \in
D_i^{\rm good}$, except for one $x_k \in D_k^{bad}$ ($i,k \in
V_0$). The next discrepancy generates $\frac{1}{2}(|V_0|^2 + |V_0|)$
subproblems, each of which contains two variables $x_{k_1} \in
D_{k_1}^{bad}$ and  $x_{k_2} \in D_{k_2}^{bad}$ ($k_1, k_2 \in V_0$),
and so on. Since we expect to obtain a very good solution 
already in the first subproblems, we will only generate subproblems up
to a certain maximum discrepancy. In our experiments the maximum
discrepancy is chosen 2 and 4 respectively.

Note that the first subproblem, corresponding to discrepancy 0,
only contains one solution, namely the one that we obtain in our
rounding procedure. By propagation of the edge constraints, all
variables $x_j \in D_j = \{0,1\}$ are instantiated automatically
to 0. Hence, the subproblem corresponding to discrepancy 0 is obsolete
in our current implementation.

\subsection{Adding discrepancy cuts}
In the case one needs to generate and solve subproblems up to a large
discrepancy, it is preferable to proof possible suboptimality of a
subproblem before entering it. Especially when the subproblems are
still relatively large. This can be done in several ways.

First, before entering a subproblem, we can identify variables which
have a subdomain of size~1, namely those with $i \in V_0$. For those
variables, one can add an extra constraint to the semidefinite
program~(\ref{eq:theta2}), being either $x_i = 1$ or $x_i = 0$. Then
the semidefinite program can be solved again, and will in general provide
a tighter bound, hopefully lower than the current lowerbound, in which
case we have proven suboptimality. However, solving the semidefinite
program each time before entering a subproblem is very time-consuming
and this method will not be very practical.

A better alternative would be to add a specific constraint, a {\em
  discrepancy cut}, that is valid for all subproblems of a given
discrepancy. Recall that $V_0 = \{i \in V \; | D_i^{\rm good} = \{1\}
\}$. Hence, all subproblems
of discrepancy $k$ consist of $k$ variables $x_i$ with $x_i = 0$, and
$|V_0| - k$ variables $x_i$ with $x_i = 1$ $(i \in V_0)$. This gives
rise to two discrepancy cuts, given discrepancy $k$: 
\begin{eqnarray}
\sum_{i \in V_0} x_i = |V_0| - k \\ \label{eq:discr_cut1}
\sum_{i \in V_0} 1 - x_i = k
\end{eqnarray}
We implemented both of them, and cut (\ref{eq:discr_cut1}) gives the
best results. Stated in terms of semidefinite
program~(\ref{eq:theta2}), the discrepancy cut looks like tr$(AX) =
|V_0| - k$, with $A_{ij} = 1$ if either $i \in V_0$ or $j \in V_0$ or
both, and $A_{ij} = 0$ otherwise. As mentioned before, solving a
semidefinite relaxation is relatively expensive, and one should make a
tradeoff between its computation time and the gain in time of not
solving the subproblem. For the instances we considered, the time
needed to solve a semidefinite program is always larger than the time
needed to solve all subdomains we would like to proof
suboptimal. However, these cuts might be helpful for larger
instances. 

\section{Computational results} \label{sc:results}
Our experiments are being done on a Pentium 1GHz processor, with
256 Mb RAM, using Windows XP. As constraint programming solver we use
the ILOG solver library, version 5.1 \cite{ilog51}. As semidefinite
programming solver, we use CSDP version 4.1 \cite{csdp}.

The first instances we consider are randomly generated weighted graphs with
$n$ vertices and $m$ edges. The vertex weights are randomly chosen
from a range of 1 up to $n$. The edge density is chosen such that the
constraint programming solver has difficulties solving
them. Namely, the more edge constraints we have, the easier it is
solved by constraint programming. The name of the instances represent
the number of vertices and the edge density, i.e. {\tt g75d015} is a
graph on 75 vertices with an approximate edge density of 0.15. For
these graphs, we have chosen to generate subproblems up to a maximum
discrepancy of 4, based upon earlier experience.

We also considered structured instances, obtained from problems
arising in coding theory~\cite{sloane}. These are unweighted graphs,
therefore we have set all weights equal to 1. For these graphs, we
generate subproblems up to a discrepancy of 2.

\begin{table}
\begin{center}
\rotatebox{90}{
\begin{tabular}{|lrr|p{1.3cm}p{1cm}p{.8cm} c p{1.1cm}p{1cm}p{1.1cm}p{1cm}|p{1cm}p{1cm}p{1.3cm}|} \hline
\multicolumn{3}{|c|}{\bf instance} & \multicolumn{8}{|c|}{\bf sdp and cp} &
\multicolumn{3}{|c|}{\bf cp alone} \\
 & & & & & & best & time & time & total & back- & & total & back-  \\
 name & $n$ & $m$ & $\vartheta$ & round & best & \hspace{.16cm}discr\hspace{.16cm} & sdp &
 subp & time & tracks & best & time & tracks \\ \hline
{\tt g50d005} & 50 & 69 & 746.00 & 746 & 746$^*$ & 0 & 0.48 & 0.00 & 0.48 & 0 & 746$^*$ & 6.38 & 160185 \\
{\tt g50d010} & 50 & 130 & 568.00 & 568 & 568$^*$ & 0 & 0.53 & 0.00 & 0.53 & 0 & 568$^*$ & 1.54 & 40878 \\
{\tt g50d015} & 50 & 191 & 512.00 & 512 & 512$^*$ & 0 & 0.71 & 0.00 & 0.71 & 0 & 512$^*$ & 0.54 & 13355 \\
{\tt g75d005} & 75 & 139 & 1472.17 & 1455 & 1466 & 1 & 0.99 & 3.47 & 4.46 & 36492 & 1077 & 100.00 & 2555879 \\
{\tt g75d010} & 75 & 280 & 1148.25 & 1122 & 1134 & 3 & 2.05 & 1.01 & 3.06 & 12964 & 1074 & 100.00 & 2299226 \\
{\tt g75d015} & 75 & 414 & 966.76 & 924 & 946 & 4 & 3.50 & 0.88 & 4.38 & 11136 & 951$^*$ & 46.63 & 1000409 \\
{\tt g100d005} & 100 & 250 & 2903.00 & 2903 & 2903$^*$ & 0 & 3.13 & 0.00 & 3.13 & 0 & 2415 & 100.01 & 1578719 \\
{\tt g100d010} & 100 & 495 & 2058.41 & 1972 & 2029 & 4 & 6.92 & 7.17 & 14.09 & 74252 & 1850 & 100.02 & 1666892 \\
{\tt g100d015} & 100 & 725 & 1704.61 & 1568 & 1608 & 4 & 17.87 & 4.28 & 22.15 & 47222 & 1644 & 100.01 & 1469221 \\
{\tt g125d005} & 125 & 367 & 3454.00 & 3454 & 3454$^*$ & 0 & 7.10 & 0.00 & 7.10 & 0 & 2656 & 100.02 & 1644993 \\
{\tt g125d010} & 125 & 761 & 2448.94 & 2208 & 2271 & 4 & 22.00 & 13.98 & 35.98 & 107610 & 1668 & 100.01 & 1643183 \\
{\tt g125d015} & 125 & 1110 & 2033.74 & 1839 & 1846 & 4 & 59.56 & 6.15 & 65.71 & 53103 & 1733 & 100.02 & 1454377 \\
{\tt g150d005} & 150 & 549 & 5043.09 & 5035 & 5035$^*$ & 0 & 15.79 & 26.89 & 42.68 & 135051 & 3090 & 100.02 & 1426234 \\
{\tt g150d010} & 150 & 1094 & 3651.38 & 3281 & 3423 & 2 & 64.83 & 13.10 & 77.93 & 78263 & 2294 & 100.02 & 1462125 \\
{\tt g150d015} & 150 & 1641 & 2935.84 & 2572 & 2572$^*$ & 0 & 181.54 & 4.56 & 186.10 & 27306 & 2224 & 182.03 & 2119929 \\
 &  &  &  &  &  &  &  &  &  &  &  &  &  \\
{\tt 1tc.64} & 64 & 192 & 20.00 & 20 & 20$^*$ & 0 & 1.05 & 0.00 & 1.05 & 0 & 19 & 324.01 & 10969580 \\
{\tt 1et.64} & 64 & 264 & 18.85 & 18 & 18$^*$ & 0 & 0.97 & 0.00 & 0.97 & 0 & 18$^*$ & 179.00 & 5763552 \\
{\tt 1tc.128} & 128 & 512 & 38.00 & 38 & 38$^*$ & 0 & 12.03 & 0.00 & 12.03 & 0 & 19 & 324.02 & 9455475 \\
{\tt 1et.128} & 128 & 672 & 29.33 & 28 & 28$^*$ & 0 & 13.48 & 0.67 & 14.15 & 1929 & 19 & 324.02 & 8562959 \\
{\tt 1dc.128} & 128 & 1471 & 16.89 & 16 & 16$^*$ & 0 & 107.71 & 0.00 & 107.71 & 0 & 13 & 324.04 & 6497027 \\
{\tt 1zc.128} & 128 & 2240 & 20.84 & 16 & 18 & 2 & 323.53 & 0.41 & 323.94 & 2094 & 18 & 324.04 & 7584769 \\
{\tt 1tc.256} & 256 & 1312 & 63.43 & 60 & 62 & 2 & 129.82 & 8.63 & 138.45 & 11588 & 13 & 324.03 & 7284451 \\
\hline
\end{tabular}
}
\end{center}
\caption{Computational results on randomly generated weighted graphs
  and structured unweighted graphs. Best found solutions that are
  proven optimal are marked with an asterisk (*).}\label{tb:results}
\end{table}

The results of our experiments are given in
Table~\ref{tb:results}. It consists of three parts: the first part
describes the instances, the next part gives the results of our
approach ({\em sdp and cp}), the last part concerns the results of a
sole constraint programming approach ({\em cp alone}).

The columns in this table represent the following.
An instance {\em name} has $n$ vertices and $m$ edges. For the part on
our approach, the value of the semidefinite relaxation is $\vartheta$,
the rounded solution of the semidefinite relaxation has value {\em
  round}, and {\em best} is the value of the best solution found. This
best solution is found in a subproblem generated during discrepancy
{\em best discr}. Note that we 
generate subproblems up to discrepancy 4 in all cases, as was
mentioned in Section~\ref{sc:framework}. The time spent on solving the
semidefinite relaxation is denoted by {\em time sdp}. The time spent
on solving all generated subproblems is denoted by {\em time
  subp}. These values together form the {\em total time}. All
times are measured in seconds. The
number of all backtracking steps made during the search in our
approach is collected in {\em backtracks}. Concerning the sole
constraint programming approach, we report the best solution found
({\em best}),  the {\em total time} spent during search, and the total
number of {\em backtracks}. Note that we have set a time limits for
the constraint programming solver, to create a fair comparison with
our approach. They are 100 seconds for {\tt g50d005} up to {\tt
  g150d010}, 190 seconds for {\tt g150d015} and 324 seconds for the
structured instances. Best found solutions that are proven to be
optimal are indicated by an asterisk (*).

For the instances in Table~\ref{tb:results} we only solved one
semidefinite relaxation per problem, namely at the root node. The
reason for this is that the time spent during the subproblem search is
less than the time spent on computing another relaxation, as reported
in the table. Therefore, we cannot gain time by adding discrepancy
cuts and computing another semidefinite relaxation.

In general, our method produces better bounds than the constraint
programming approach alone. In many cases, the rounded solution of the
semidefininite relaxation is already optimal. However, note that there
are two instances for which the constraint programming approach gives
better bounds. Note also that the structured instances are handled
quite well by our approach, while the constraint programming approach
produces very low quality solutions for the larger instances.

\section{Related literature}\label{sc:literature}
Since the stable set problem is NP-hard, no complete (or exact)
algorithm is known that solves the stable set problem in polynomial
time. Many other techniques have been proposed, including
approximation algorithms, heuristics, or branch and bound structured
methods. A survey of different formulations, complete methods and
heuristics for the maximum clique problem\footnote{The maximum
  weighted stable set problem of a graph is equivalent to the maximum
  weighted clique problem of its complement graph.}
is given by Pardalos and Xue~\cite{pardalos_xue} and, more recently,
by Bomze et al.~\cite{bomze}. 

Although semidefinite programs can be solved in polynomial time
theoretically, it lasted until a few years ago until fast solvers for
this purpose were implemented. Until then, application inside a branch
and bound framework was unrealistic. Still, solving a semidefinite
program takes relatively much time, compared to solving a linear
program. However, since semidefinite programming solvers are getting
faster, semidefinite relaxations become a serious candidate to be used
within a branch and bound framework, see for instance the paper by
Karisch et al. \cite{karisch2000}. 

A large number of references to papers concerning semidefinite
programming are on the web pages of Helmberg~\cite{sdp_helmberg} and
Alizadeh~\cite{sdp_alizadeh}. A general introduction on semidefinite
programming applied to combinatorial optimization is given by Goemans
and Rendl~\cite{goemans_rendl}.

Another area that made semidefinite programming useful in practice is
that of approximation algorithms. In this field one tries to give
a performance guarantee for an algorithm on a particular problem. In
particular, the paper \cite{GW_maxcut} by Goemans and Williamson uses
a semidefinite relaxation and randomized rounding to prove such a
performance guarantee for the maximum cut problem of a graph and
satisfiability problems. Following this result numerous papers
appeared, also concerning the approximation of satisfiability
problems, including \cite{zwick4sat} and \cite{lau2002}.

The solution structure of the current work, namely problem
decomposition by branching on promising subdomains, is similar to the
method described in~\cite{milano_hoeve02}. In that paper a linear
relaxation is used to identify promising values. Moreover, by
exploiting the discrepancy structure of the method combined with
reduced costs, suboptimality of subproblems can be proved very fast.

Another hybrid approach using linear programming and constraint
programming has been investigated by Ajili et al.~\cite{ajili01lp} and
El Sakkout et al.~\cite{sakkout_probe}.
A subset of constraints is relaxed as a linear program in such a
way that its solution is always integral. The solution to the
relaxation serves as a suggestion (a `probe') for solving the complete
program using a constraint programming solver. A probe is used to
detect infeasibility, to remove inconsistent domain values and to
guide the search. During search, many probing steps are being
made. This results in a tight cooperation of the linear programming
and constraint programming solver.

\section{Conclusions and further research} \label{sc:conclusions}
We introduced a method that combines semidefinite programming and
constraint programming to solve the stable set problem. Our
experiments show that constraint programming can indeed benefit
greatly from semidefinite programming. On instances that were very
difficult to handle for a constraint programming solver, our hybrid
method obtained very good results.

The discrepancy cuts we proposed to strengthen the semidefinite
relaxation could not be applied efficiently to the instances we
considered. However, for larger instances they could be helpful.

Further research in this direction would for instance be to obtain a
filtering mechanism similar to the cost-based domain filtering for
linear relaxations \cite{NostroCP99}. In~\cite{helmberg}, Helmberg
describes such a procedure, called variable fixing, for semidefinite
relaxations. It would be interesting to see how his method can be
applied in a constraint programming framework. 

Also, one could consider a different way of selecting promising
values from the solution of the semidefinite relaxation. A strategy
that incorporates randomized rounding possibly yields better
results. This thought is motivated by the use of randomized rounding
of semidefinite relaxations in approximation algorithms, as discussed
in Section~\ref{sc:literature}.

Finally, this work has much in common with our previous work
\cite{milano_hoeve02}. The underlying general principle of decomposing
a problem into promising subproblems according to a certain heuristic
is currently under research.

\section*{Acknowledgements}
Many thanks to Michela Milano for fruitful discussion and helpful
comments on this paper.


\end{document}